\documentclass{article}

\usepackage{amssymb, latexsym}

\usepackage{graphicx}

\usepackage{amsmath}

\newtheorem{theorem}{Theorem}[section]

\newtheorem{problem}[theorem]{Problem}

\newtheorem{remark}[theorem]{Remark}

\date{ }

\begin{document}
\title{\Large\bf  From the Boltzmann $H$-theorem to Perelman's $W$-entropy formula for the Ricci flow\footnote{Appeared in Proceedings on Sino-Japan Conference of Young Mathematicians on Emerging Topics on Differential Equations and their Applications, Chern Institute of Mathematics, December 5-9, 2011, Tianjin, Nankai University. Edited by Hua Chen, Yiming Long and Y. Nishiura, World Scientific Publishing, pp. 68-84.}}
\author{Xiang-Dong Li\thanks{Research supported by NSFC No. 10971032, Key Laboratory RCSDS, CAS, No. 2008DP173182, and a
Hundred Talents Project of AMSS,
CAS.}
}

\maketitle

\medskip

\begin{minipage}{120mm}
{\bf Abstract}. In 1870s, L. Boltzmann proved the famous
$H$-theorem for the Boltzmann equation in the kinetic theory of
gas and gave the statistical
interpretation of the thermodynamic entropy.  In 2002, G. Perelman introduced the notion of $W$-entropy and proved the $W$-entropy formula for the Ricci
flow. This plays a crucial role in the
proof of the no local collapsing theorem and in the final resolution
of the Poincar\'e conjecture and Thurston's geometrization
conjecture. In our previous paper \cite{Li11a}, the author gave a probabilistic interpretation of the $W$-entropy using
the Boltzmann-Shannon-Nash entropy. In this paper, we  make some further efforts for a better understanding of the mysterious $W$-entropy by comparing the $H$-theorem for the Boltzmann equation and the Perelman $W$-entropy
formula for the Ricci flow. We also suggest a way to construct the
``density of states'' measure for which the Boltzmann $H$-entropy is exactly the
$W$-entropy for the Ricci flow.
\end{minipage}

\section{Introduction}

In a seminal paper \cite{P1}, Perelman introduced the $W$-entropy
for the Ricci flow and the associated conjugate heat equation and proved its monotonicity  on
compact manifolds. More precisely, let
$M$ be an $n$-dimensional compact manifold with a family of
Riemannian metrics $\{g(t): t\in [0, T]\}$,  $f\in C^\infty(M\times
[0, T], \mathbb{R})$, where $T>0$ is fixed. Suppose that
 $(g(t), f(t), \tau(t), t\in [0, T])$ is a solution to the evolution
equations
\begin{eqnarray}
\partial_t g&=&-2Ric,\label{RF}\\
\partial_t f&=&-\Delta f+|\nabla
f|^2-R+\frac{n} {2\tau},\label{r-c}\\
\partial_t \tau&=&-1,\label{tau}
\end{eqnarray}
Following \cite{P1}, the $W$-entropy functional associated to
$(\ref{RF})$ and $(\ref{r-c})$ is defined by
\begin{eqnarray}
W(g, f, \tau)=\int_M \left[\tau(R+|\nabla
f|^2)+f-n\right]\frac{e^{-f}} {
(4\pi\tau)^{n/2}}dv,\label{entropy-1}
 \end{eqnarray}
where $v$ denotes the volume measure on $(M, g(\tau))$. By
\cite{P1}, the following entropy formula holds
\begin{eqnarray}
\frac{d} {dt}W(g, f, \tau)=2 \int_M \tau\left|Ric+\nabla^2
f-\frac{g} {2\tau}\right|^2\frac{e^{-f}}{(4\pi
\tau)^{n/2}}dv.\label{ent-p}
\end{eqnarray}
As a consequence of the Perelman $W$-entropy formula
$(\ref{ent-p})$, the $W$-entropy is monotone increasing in $\tau$
and the monotonicity is strict except that $(M, g(\tau), f(\tau))$
is a shrinking Ricci soliton, i.e.,
\begin{eqnarray}
Ric+\nabla^2 f={g\over 2\tau}.\label{SRS}
\end{eqnarray}
This result plays a crucial role in the proof of the no local
collapsing theorem, which is equivalent to the long time standing
Hamilton's Little Loop Conjecture \cite{H2}. As pointed out by
Perelman \cite{P1}, ``this removes the major stumbling block in
Hamilton's approach to geometrization'' (for the resolution of the
Poincar\'e conjecture and Thurston's geometrization conjecture). See
also \cite{CZ, CLN, KL, MT}.

Since Perelman's paper \cite{P1} was posted on Arxiv in 2002, many
people have established the $W$-entropy formula for various
geometric evolution equations \cite{CKL, Li11a, Li11b, N1,N2, Ec,
LNVV, KN}. However, it seems that there is no reference in the
literature which explains clearly how the mysterious $W$-entropy was
introduced for the Ricci flow and what are the hidden insights for
Perelman to introduce the $W$-entropy for the Ricci flow. In our previous paper \cite{Li11a}, the author gave a probabilistic interpretation of the $W$-entropy using the
the Boltzmann-Shannon-Nash entropy. The purpose of this paper
is try to make some further efforts for a better understanding of
the mysterious $W$-entropy by comparing the $H$-theorem for the Boltzmann equation and the Perelman $W$-entropy
formula for the Ricci flow. We also suggest a way to construct the
``density of states'' measure for which the Boltzmann $H$-entropy is exactly the
$W$-entropy for the Ricci flow.

\section{History and Boltzmann entropy formula}

\subsection{Thermodynamic entropy}

The notion of entropy was introduced by R. Clausius in 1865 in his
study of the Carnot cycle in the thermodynamics \cite{Cl1, Cl2}. He
proved that, in any reversible thermodynamic process, ${\delta
Q\over T}$ is an exact form, where $\delta Q$ denotes the total
change of the heat along the reversible thermodynamic process, and
$T$ is the temperature. This leads Clausius to define the
thermodynamic entropy $S$ as a state function of the thermodynamic
system satisfying
\begin{eqnarray*}
dS={\delta Q\over T}.
\end{eqnarray*}
Thus, if the reversible thermodynamic process changes its states from $A$ to $B$, then
\begin{eqnarray*}
S(B)-S(A)=\int_A^B {\delta Q\over T}.
\end{eqnarray*}
Moreover, the value of $S(B)-S(A)$ is independent of the choice of
any reversible process from $A$ to $B$. If the thermodynamic process
is irreversible, the second law of thermodynamics says that: The
thermodynamic entropy $S$ of an isolated thermodynamic system always
increases, i.e.,
$$\Delta S\geq 0.$$
In this sense, the thermodynamic entropy $S$ is an important
quantity (as a function of the state)  to characterize the
irreversibility of the thermodynamic processes.

%By the basis formula of the thermodynamics
%\begin{eqnarray*}
%\delta Q=dU+pdV,
%\end{eqnarray*}
%where $U$ is the inner energy, $p$ is the thermodynamic pressure, $V$ the volume of ideal gas, and $T$ the temperature. This leads to
%\begin{eqnarray*}
%dS={1\over T}dU+{1\over T}pdV,
%\end{eqnarray*}

%Let $S$ be the thermodynamic entropy, $U$ the inner energy of a thermodynamic system, $T$ the temperature, and $V$ the volume of gas. The Helmholtz free energy is defined by
%\begin{eqnarray*}
%F(T, V)=U(S, V)-TS.
%\end{eqnarray*}
%from which one has
%\begin{eqnarray*}
%S=-{\partial F\over \partial T}.
%\end{eqnarray*}

\subsection{Boltzmann equation and $H$-theorem}

In 1872, Boltzmann published an important paper \cite{Btz1}, which
contained two celebrated results nowadays known as the Boltzmann
equation and the Boltzmann $H$-theorem. More precisely, let
$f(x, v, t)$ be the probabilistic distribution of the ideal gas at
time $t$ in the phase space $\mathbb{R}^6$ of position $x$ and
velocity $v$, then $f$ satisfies the Boltzmann equation (see e.g.
\cite{Btz1, Btz3, Ev, CIP, C1, C2})
\begin{eqnarray*}
\partial_t f+v\cdot \nabla_x f=Q(f, f)
\end{eqnarray*}
where
\begin{eqnarray*}
Q(f, f)= \int_{\mathbb{R}^3}\int_{S^2}[f(x, v')f(x, v'_*)-f(x,
v)f(x, v_*)]B(v-v_*, \theta)dv_* dS(u),
\end{eqnarray*}
where $dS$ denotes the surface measure on $S^2$, $v', v'_*$ are
defined in terms of $v, v_*, u$ by
\begin{eqnarray*}
v'=v-[(v-v_*)\cdot u]u,  \ \ \ v_*'=v_*+[(v-v_*)\cdot u]u,
\end{eqnarray*}
and $B: \mathbb{R}^3\times S^2\rightarrow (0, \infty)$ is the
collision kernel and is assumed to be rotationally invariant, i.e.,
$B(z, u)=B(|z|, |z\cdot u|)$, $\forall z\in \mathbb{R}^3$ and $u\in
S^2$. Let
\begin{eqnarray}
H(t)=\int_{\mathbb{R}^6}f(x, v, t)\log f(x, v, t)dxdv \label{BoH}
\end{eqnarray}
be the Boltzmann $H$-functional of the probability distribution
$f(x, v, t)dxdv$. The Boltzmann $H$-Theorem states that, if $f(x, v,
t)$ is the solution of the Boltzmann equation, which is
``sufficiently well-behaved'' (in the sense of Villani \cite{V1},
for its precise meaning, see \cite{CIP}), then the following formula
holds (see e.g. \cite{Btz1, Btz3, Ev, CIP, C1, C2})
\begin{eqnarray}
{dH\over dt}={1\over 4}
\int_{S^2}\int_{\mathbb{R}^3}(f'f'_*-ff_*)[\log (ff_*)-\log
(f'f'_*)]B(v-v_*, \theta)dv_* dS(u),\label{BDF}
\end{eqnarray}
where $f=f(\cdot, v, \cdot)$, $f'=f(\cdot, v', \cdot)$,
$f_{*}=f(\cdot, v_*, \cdot)$, $f'_{*}=f(\cdot, v'_{*}, \cdot)$. By
the Boltzmann formula for the $H$-functional $(\ref{BDF})$ and using
the elementary inequality
\begin{eqnarray*}
(x-y)(\log y-\log x)\leq 0, \ \ \ \forall x, y\in \mathbb{R}^+,
\end{eqnarray*} one can conclude that
$H$ is always decreasing in time, i.e.,
\begin{eqnarray}
{dH\over dt}\leq 0, \ \ \ \forall t>0,\label{HTh}
\end{eqnarray}
and the equality holds in $(\ref{HTh})$ if and only if
\begin{eqnarray}
f'f'_{*}=ff_*,\ \ \ \forall v, v_*\in \mathbb{R}^3, u\in
S^2,\label{EBE}
\end{eqnarray}
 which
implies that $f$ is a local Maxwell distribution, i.e.,
\begin{eqnarray}
f(x, v, t)=n(x, t)\left({m\over 2\pi k
T(x, t)}\right)^{3/2}\exp\left(-{|v-\bar v(x, t)|^2\over 2kmT(x, t)}\right),\label{Max}
\end{eqnarray}
where the parameters $m$ is the mass of particle, $n(x, t)\in \mathbb{R}^3$ is the particle density at $(x, t)$,
$\bar v(x, t)\in \mathbb{R}^3$ and $T(x, t)> 0$ are
the mean velocity and the local temperature at $(x, t)$.

\subsection{Boltzmann entropy formula}

In 1877, to better explain his $H$-theorem, Boltzmann
\cite{Btz2} introduced the statistical interpretation of the
thermodynamic entropy $S$ by the formula
\begin{eqnarray}
S= k \log W,\label{SkW}
\end{eqnarray}
where $k$ is a constant, which is nowadays called the Boltzmann
constant (for simplicity, we choose $k=1$ throughout this paper),
and $W$ is the German word ``Wahrscheinlichkeit'' of the English
word ``probability'', which is identified with the number of
possible microstates corresponding to the macroscopic state of a
given thermodynamic system, i.e., the number of (unobservable)
``way'' in the (observable) thermodynamic state of a system can be
realized by assigning different positions and momenta to various
molecules of the ideal gas. The Boltzmann entropy formula
$(\ref{SkW})$ indicates the logarithmic connection between Clausius'
thermodynamic entropy $S$ and the number $W$ of the most probable
microstates consistent with the given macrostate.

We now explain the Boltzmann entropy formula $(\ref{SkW})$ by
studying the following question: Take $N$ identical particles, which
will be distributed over $k$ boxes, and let $p_1, \ldots, p_k$ be
some rational numbers in $[0, 1]$ such that
$\sum\limits_{i=1}^kp_i=1$. We can regard $p_j$ as the probability
for each particle to be distributed into the $j$-th box. Let $W$ be
the number of configurations such that $N_j$ particles are
distributed in the $j$-th box, $j=1, \ldots, k$, and $P$ the
corresponding probability. Then
\begin{eqnarray*}
W={N!\over N_1!\ldots N_k!},
\end{eqnarray*}
and
\begin{eqnarray*}
P={N!\over N_1!\ldots N_k!}p_1^{N_1}\ldots p_k^{N_k}.
\end{eqnarray*}
Assuming that $N_j>>1$, $j=1, \ldots, k$. By the Stirling formula $\log N!\equiv (N+1/2)\log N+N+{1\over 2}\log (2\pi)+O(1/N)$, we have
\begin{eqnarray*}
\log P\simeq \left(N+{1\over 2}\right)\log N-\sum\limits_{i=1}^k \left(N_j+{1\over 2}\right)\log N_j-\sum\limits_{j=1}^k N_j \log p_j-{k-1\over 2}\log(2\pi).
\end{eqnarray*}
Using the method of the Lagrangian multiplier, we can prove that the maximum value of $P$ subject to the constraint $\sum\limits_{j=1}^kN_j=N$ is achieved at $N_j=p_jN$ and we have
\begin{eqnarray*}
\log P_{\rm max}\simeq -{1\over 2}\log N-\sum\limits_{i=1}\left(2p_j N+{1\over 2}\right)\log p_j-{k-1\over 2}\log (2\pi).
\end{eqnarray*}
Therefore
\begin{eqnarray*}
\lim\limits_{N\rightarrow \infty}{\log P_{\rm max}\over N}=-2\sum\limits_{i=1}^k p_i\log p_i.
\end{eqnarray*}
Let $W_{\rm max}$ be the corresponding number of the configurations
such that $P$ achieves its maximum value $P_{\rm max}$. Then (cf.
also \cite{V1})
\begin{eqnarray}
S:=\lim\limits_{N\rightarrow \infty}{\log W_{\rm max}\over
N}=-\sum\limits_{i=1}^k p_i\log p_i.\label{S}
\end{eqnarray}

By comparing  the Boltzmann entropy formula $(\ref{S})$ with the
$H$-functional $(\ref{BoH})$ introduced by Boltzmann in the study of
the kinetic theory of gas, one can see that the $H$-functional is
nothing else but the minus entropy of the probability distribution
$f(\cdot, \cdot, t)$ in the phase space at time $t$.

\subsection{Shannon and Nash entropy}
In 1948, C. Shannon \cite{Shan} introduced the notion of entropy
into the theory of information. More precisely,
$$
H=-\sum\limits_{i=1}^k p_i \log p_i.$$
According to Shannon \cite{Shan}, ``the form of $H$ will be recognized as that of entropy as defined in certain formulations of statistical mechanics where $p_i$ is the probability of a system being in cell $i$ of its phase space... $H$ is then, for example, the $H$-in Boltzmann's famous $H$-theorem.''

On the other hand, J. Nash \cite{Na} used the $H$-entropy to study
the continuity of the parabolic and elliptic PDEs. In the
literature, the Boltzmann $H$-entropy is sometime called the Shannon
entropy or the Nash entropy. In this paper, we call $H$ the
Boltzmann-Shannon-Nash entropy, or for simplicity the Boltzmann
entropy.

\subsection{The maximum entropy principle and the central limit theorem}

According to Boltzmann's entropy formula, the entropy $S$ of a given
macrostate achieves its maximum value at the most probable
microstates subject to the constraint $\sum\limits_{j=1}^k N_j=N$.

The following theorem, which can be proved by the method of Lagrangian multiplier, indicates a deep connection between the central limit theorem and the maximum entropy principle, in view of which, one can interpret the central limit theorem as a consequence of the maximum entropy principle.

\begin{theorem} Let $\mathcal{D}=\{\mu=fdx\in \mathbb{P}(\mathbb{R}): \int_M xd\mu(x)=0, \ \int_M x^2d\mu(x)=1\}$. Then
\begin{eqnarray*}
\gamma=\arg \min\{H(\mu): \mu\in \mathcal{D}\}.
\end{eqnarray*}
\end{theorem}

\subsection{Canonical ensemble and Boltzmann's entropy formula}

%The canonical ensemble in statistical mechanics is a statistical ensemble representing a probability distribution of microscopic states of the system.
%From {\it Wikipedia}, we can find the following definition \footnote{Eric Carlen pointed to the author that this definition can apply to any of the micro-canonical, canonical
%and also grand-canonical  distributions. The canonical distribution represents the probability distribution of the microstates of a system in equilibrium and in contact with a heat-bath at inverse temperature $\beta$. The microcaninical is for an isolated system, and the grand-canonical for a system in equilibrium and in contact with a heat bath and a particle reservoir.}
%: ``The
%canonical ensemble in statistical mechanics is a statistical
%ensemble representing a probability distribution of microscopic
%states of the system.''
Heuristically, the canonical ensemble is a measurable  configuration space $(\Omega,
\mathcal{F})$ on which ``there exist'' {\it a priori  measure}  ${\it D}\omega$ and a Hamiltonian function
$E: \Omega\rightarrow [0, +\infty]$ , such that the distribution of the particles on $(\Omega, \mathcal{F})$ is given by a probability measure $\mathbb{P}$ on $(\Omega, \mathcal{F})$
which has the following formal  expression

$$
d\mathbb{P}(\omega)={1\over Z_\beta}e^{-\beta E(\omega)}{\it D}\omega,$$
where
$$Z_\beta=\int_\Omega e^{-\beta E(\omega)}{\it D}\omega,$$
is called the partition function. In the literature, the probability measure $\mathbb{P}$ is called the Gibbs measure. We use $(\Omega, \mathcal{F}, E, \mathbb{P})$ to denote the canonical ensemble $(\Omega, \mathcal{F})$ with the Hamitonian function $E$ and with the Gibbs measure $\mathbb{P}$. In the kinetic theory of ideal gas, $\Omega$ is the phase space of the positions and the velocities of the particles, and $\mathbb{P}$ is also called the Maxwell-Bolztmann distribution.

To introduce the Maxwell-Boltzmann distribution, let us consider the
following problem: Let us assume that $N$ particles of ideal gas are
distributed into $n$ rooms $R_1, \ldots, R_n$, each room $R_i$ has
$g_i$ boxes, $i=1, \ldots, n$. Suppose that if a particle is in the
room $R_i$, it has energy $E_i$, $i=1, \ldots, n$. Then the total
number of ways that these $N$ particles are distributed in to boxes
of these $n$ rooms, for which $N_i$ particles are in the room $R_i$
and with energy $E_i$, is given by
$$
W={N! \over \prod_{i=1}^n N_i!} \prod_{i=1}^n g_i^{N_i}.$$

The most probable distribution is the one for which $W$ achieves its
maximum value under two constraint conditions:

$$N=\sum\limits_{i=1}^n N_i,\ \ \ \ E=\sum\limits_{i=1}^n N_i E_i.$$

By the method of Lagrange  multiplier, and using the Stirling
asymptotic formula, one can prove that, in the case $N>>1$ and
$N_i>>1$, the most probable distribution satisfying the above two
constraint conditions is given by
$$
N_i=A g_i e^{-\beta E_i},,  \ \ \ \ i=1, \ldots, n.$$ where $A$ and
$\beta$ are two constants. By the theory of thermodynamic, it is known that $\beta$ is related
to the absolute temperature  $T$ via
$$\beta={1\over kT},$$
where $k$ is a constant, called the Boltzmann constant. To determine
$A$, it is useful to introduce the partition function
$$
Z_\beta=\sum\limits_{i=1}^n g_i e^{-\beta E_i}.$$ Then
$$A={N\over Z_\beta}.$$

Hence, the most probable distribution, i.e., the Maxwell-Boltzmann
distribution, is given by
$$
{N_i\over N}={1\over Z_\beta}g_i e^{-\beta E_i}, \ \ \ \ i=1,
\ldots, n.$$ By Boltzmann's entropy formula $(\ref{SkW})$ in
Statistical Mechanics, the thermodynamic entropy of the macrostate
system is given by
$$S=k \log W.$$
Using the explicit expression of the Maxwell-Boltzmann distribution,
one can prove that
$$
S=k N\log Z_\beta+k\beta E.$$ Let
$$\langle E\rangle=\sum\limits_{i=1}^n {n_i E_i\over N}.$$  Then we obtain the Boltzmann statistical mechanics
interpretation of the $S$-entropy:
$$S=\lim\limits_{N\rightarrow \infty}{\log W\over N}=\log Z_\beta+\beta \langle E\rangle.$$

In the continuous case, let  $\mathbb{P}$ be the Gibbs measure on a canonical ensemble $(\Omega,
\mathcal{F})$ with Hamiltonian energy function $E:\Omega\rightarrow \mathbb{R}^+$. Then the Maxwell-Boltzmann distribution is given by the following probability
measure on $(\mathbb{R}^+, \mathcal{B}(\mathbb{R}^+))$
\begin{eqnarray*}
dP(E)={e^{-\beta E}\over Z_\beta}g(E)dE,
\end{eqnarray*}
where $\beta\in \mathbb{R}$ is the Boltzmann constant, $g(E)dE$
denotes the ``density of states'' measure, whose physical meaning is
the number of microstates with energy levels in the range $[E,
E+dE]$, and
\begin{eqnarray*}
Z_\beta=\int_{\mathbb{R}^+} e^{-\beta E}g(E)dE
\end{eqnarray*}
is the partition function. In this case, one can define the
temperature of the macrostate system by
\begin{eqnarray*}
T={1\over \beta},
\end{eqnarray*}
and define the Helmholtz  free energy by
\begin{eqnarray*}
F=-{1\over \beta}\log Z_\beta.
\end{eqnarray*}
One can verify that, the average of the energy with
respect to the Maxwell-Boltzmann distribution satisfies
$$\langle E\rangle=-{\partial\over \partial \beta}\log Z_\beta,$$
and the following formula holds for the Boltzmann entropy $S$ of the
Maxwell-Boltzmann distribution:
$$S=\log Z_\beta-\beta {\partial \over \partial \beta}\log Z_\beta.$$
Moreover, the fluctuation of the energy with respect to the
Maxwell-Boltzmann distribution is given by
$$
\sigma:=\langle (E-\langle E\rangle)^2\rangle ={\partial^2 \over
\partial \beta^2} \log Z_\beta,$$
and the derivative of the entropy with respect to $\beta$ satisfies
\begin{eqnarray*}
{\partial S\over \partial \beta}=-\beta{\partial^2\over \partial
\beta^2}\log Z_\beta=-\beta \sigma.
\end{eqnarray*}
By the definition formula $\beta={1\over T}$, we can verify the following formulas

\begin{eqnarray*}
\langle E\rangle&=&T^2{\partial\over \partial T}\log Z_\beta,\ \ F=-T\log Z_\beta,\\
S&=&{\partial \over \partial T}(T\log  Z_\beta)=-{\partial F\over \partial
T},
\end{eqnarray*}
and
$$
{\partial S\over \partial T}={\sigma\over T^3}.$$

\section{Perelman's interpretation of the $W$-entropy for Ricci flow}

In Section $5$ in \cite{P1}, Perelman gave a heuristic
interpretation for the $W$-entropy using statistical mechanics. Let
$(M, g(\tau))$ be a family of closed Riemannian manifolds,
$dm(\tau)=(4\pi\tau)^{-n/2}e^{-f(\tau)}dv_{g(\tau)}$ a probability
measure on $(M, g(\tau))$, where $g(\tau)$ satisfies the Perelman's
backward Ricci flow equation ${\partial\over \partial
\tau}g_{ij}=R_{ij}+\nabla_{i}\nabla_{j}f$, and $f(\tau)$ satisfies
the heat equation $\partial_\tau f=\Delta f-|\nabla f|^2+R-{n\over
2\tau}$. Assume that for some canonical ensemble $(\Omega,
\mathcal{F}, E, \mathbb{P})$ on which there is a ``density of state
measure'' $g(E)dE$ such that the partition function
$Z=\int_{\mathbb{R}^+} e^{-\beta E}g(E)dE$ is given by
\begin{eqnarray}
\log Z=\int_M (-f+{n\over 2})dm,\label{logZ}
\end{eqnarray}
where $\beta={1\over \tau}$, and $\tau$ is regarded as the
temperature. Then, using the formulas in statistical mechanics, we
can obtain
\begin{eqnarray*}
\langle E\rangle&=&=-\tau^2\int_M (R+|\nabla f|^2-{n\over
2\tau})dm,\\
S&=&-\int_M (\tau(R+|\nabla f|^2)+f-n)dm,\\
\sigma&=&2\tau^4\int_M |R_{ij}+\nabla_i\nabla_j f-{1\over
2\tau}g_{ij}|^2dm.
\end{eqnarray*}
Alternatively, we can consider the evolution equations $(\ref{RF})$
and $(\ref{r-c})$ by replacing the $t$-derivatives by  minus of the
$\tau$-derivatives, where $\tau=T-t$. This implies that
\begin{eqnarray}
W=-S,\label{W1}
\end{eqnarray}
and
\begin{eqnarray}
{dW\over dt}=2 \int_M \tau\left|Ric+\nabla^2 f-{g\over
2\tau}\right|^2dm.\label{W2}
\end{eqnarray}

The above statistical mechanics interpretation $(\ref{W1})$ of the
$W$-entropy and the derivation of the derivative formula
$(\ref{W2})$ are heuristic. The problem whether there is a canonical
ensemble $(\Omega, \mathcal{F}, E, \mathbb{P})$, equivalently,
whether there is a ``density of states'' measure $g(E)dE$, such that
the partition function $Z=\int_{\mathbb{R}^+} e^{-\beta E}g(E)dE$
satisfies Perelman's required condition $\log Z=\int_M ({n\over
2}-f)dm$ remains open. It is a very interesting to know whether one
can use some ideas of the quantum field theory to construct a
``density of states'' measure on some natural canonical ensemble
such that the Boltzmann entropy of the corresponding
Maxwell-Boltzmann distribution is exactly given by $S=-W$. See
Section $5$ below.

\section{A probabilistic interpretation of the $W$-entropy for Ricci flow}

In \cite{Li11a}, the author gave a probabilistic interpretation of
Perelman's $W$-entropy for the Ricci flow. We now present this
interpretation. For this, observe that
\begin{eqnarray*}
\int_M ({n\over 2}-f)dm={n\over 2}(1+\log (4\pi \tau))-H(m),
\end{eqnarray*}
where
\begin{eqnarray*}
H(m)=-\int_M u\log udv=\int_M (f+{n\over 2}\log(4\pi
\tau)){e^{-f}\over (4\pi\tau)^{n/2}}dv
\end{eqnarray*}
is the Boltzmann-Nash-Shannon entropy of the probability measure
$dm=udv$, where $u={e^{-f}\over (4\pi\tau)^{n/2}}$. On the other
hand, let
$$d\gamma_n(x)={e^{-{\|x\|^2\over 4\tau}}\over (4\pi \tau)^{n/2}}dx$$
be the Gaussian measure on $\mathbb{R}^n$. Then it is well-known
that the Boltzmann-Nash-Shannon entropy of $\gamma_n$ is given by
\begin{eqnarray*}
H(\gamma_n)={n\over 2}(1+\log (4\pi \tau)).
\end{eqnarray*}
Hence, $\log Z=H(\gamma_n)-H(m)$ is the difference of the
Boltzmann-Nash-Shannon entropy of the Gaussian measure $\gamma_n$ on
$\mathbb{R}^n$ and the Boltzmann-Nash-Shannon entropy of the heat
kernel measure $dm=u(\tau)dv_{g(\tau)}$ on $(M, g(\tau))$. In view
of this, we have the following probabilistic interpretation of the
$W$-entropy for the Ricci flow
\begin{eqnarray*}
W={d\over d\tau}(\tau(H(\gamma_n)-H(u)),
\end{eqnarray*}
and
\begin{eqnarray*}
{dW\over d\tau}%={\sigma\over \tau^3}
=2\tau\int_M
|R_{ij}+\nabla_i\nabla_j f-{1\over 2\tau}g_{ij}|^2dm.
\end{eqnarray*}

\begin{remark}{\rm During the preparation of this paper, Songzi Li pointed out to the author that the quantity $F:=-\tau\log Z$ appeared in the definition formula of
the $W$-entropy  $W=-{\partial F\over \partial \tau}$ is actually
the Helmholtz free energy function on the canonical ensemble
$(\Omega, \mathcal{F}, E, \mathbb{P})$, where $\tau={1\over \beta}$
is regarded as the temperature. According to Otto \cite{Ot1}, if one
chooses $F(u)=\int_M u\log u dv+\int_M Vudv$ as the Helmholtz
 free energy of a certain thermodynamic system, where $V\in C^2(M)$, then the gradient flow of $F$
in the Wasserstein space equipped with a suitable infinite
dimensional Riemannian metric is the Fokker-Planck equation
$\partial_t u=\Delta u+\nabla\cdot(u\nabla V)$. In view of this, the
$W$-entropy for the Ricci flow can be regarded as the thermodynamic
entropy corresponding to the Helmholtz free energy function
$F(u)=\tau (H(\gamma_n)-H(u))$, which is the difference of the free
energy of a certain thermodynamic system on $(M, g(\tau))$ and the
free energy of the Maxwell thermodynamic equilibrium state on
$\mathbb{R}^n$. This explanation is indeed equivalent to the above
mentioned Perelman's heuristic interpretation, and is very close to
the kinetic interpretation given by the author in \cite{Li11b} for
the $W$-entropy for the Fokker-Planck equation on complete
Riemannian manifolds.}

\end{remark}

\section{Comparison between Boltzmann's $H$-theorem and
Perelman's entropy formula}

In the above sections, we have discussed, following the idea of
Perelman, the relationship between the Boltzmann $H$-entropy
and the Perelman $W$-entropy. In this section we would like to
compare the Boltzmann $H$-theorem and the Perelman $W$-entropy
formula.

The Boltzmann $H$-entropy formula $(\ref{BDF})$ and the Perelman
$W$-entropy formula $(\ref{ent-p})$ are in the same spirit in the
following three points.

(i) The Boltzmann $H$-entropy formula $(\ref{BDF})$ gives the time
derivative formula of the Boltzmann $H$-entropy along the solutions
of the Boltzmann equation, and  the Perelman $W$-entropy formula
$(\ref{ent-p})$ gives the time derivative formula of the Perelman
$W$-entropy along the solutions of the Ricci flow equation and the
conjugate heat equation.

(ii) From the Boltzmann $H$-entropy formula $(\ref{BDF})$ and
the Perelman $W$-entropy formula $(\ref{ent-p})$, we can easily
derive the monotonicity of the Boltzmann $H$-entropy along the
Boltzmann equation and the monotonicity of the Perelman $W$-entropy
along the Ricci flow and the conjugate heat equation.

(iii) We can derive the equation of the equilibrium state for the
Boltzmann $H$-entropy and for the Perelman $W$-entropy. The
equilibrium state of the Boltzmann $H$-entropy satisfies the
equation $(\ref{EBE})$, and the equilibrium state of the Perelman
$W$-entropy is the shrinking Ricci solitons which satisfy the
equation $(\ref{SRS})$.

In the case of the Boltzmann equation, the problem of the
convergence rate of the solutions of the Boltzmann equation towards
the equilibrium state (i.e., the Maxwell distribution) is the
well-known Cercignani conjecture in the study of the Boltzmann
equation. See \cite{C3, CC1, CC2, DMV, V3}. From the above
comparison, we can raise the following interesting problem for the
further study of the Perelman's $W$-entropy.

\begin{problem} What is the longtime behavior of the evolution equations $(\ref{RF})$
and $(\ref{r-c})$?  If the solution $(g(t), f(t))$ of $(\ref{RF})$
and $(\ref{r-c})$ exists on $[0, \infty)$, what is the convergence
behavior of $(g(t), f(t))$ towards its limit as $t\rightarrow
\infty$?
\end{problem}

In \cite{P1}, Perelman raised the following question

\begin{problem} If the flow is defined for all sufficiently large $\tau$ (that is, we have an ancient solution to the Ricci flow, in Hamilton's terminology), we may be interested in the behavior of the entropy $S$ as $\tau\rightarrow \infty$. A natural question is whether we have a gradient shrinking soliton whenever $S$ stays bounded.
\end{problem}
\section{An open problem: the ``density of states'' measure for the Ricci flow}

Following Perelman's heuristic interpretation of the $W$-entropy for
the Ricci flow,  we would like to discuss the following open
problem.

\begin{problem}{\rm Can one prove the existence of a certain canonical
ensemble $(\Omega, \mathcal{F})$ with {\it a priori} measure ${\it D}\omega$ and a Hamiltonian
function $E: \Omega\rightarrow \mathbb{R}^+$ such that the
partition function of the corresponding Gibbs
measure $\mathbb{P}$ on $(\Omega, \mathcal{F})$ satisfies
Perelman's condition
$$
\log \int_{\Omega}e^{-\beta E(\omega)}{\it D}\omega=\int_M \left({n\over 2}-f\right)dm?
$$
Equivalently,
$$
\log \int_{\Omega}e^{-\beta E(\omega)}{\it D}\omega=H_n(u, t),
$$
where
$$
H_n(u, t):={n\over 2}(1+\log(4\pi\tau))-\int_M u\log u dv.$$
Moreover, if such a canonical ensemble  $(\Omega, \mathcal{F}, \mathbb{P})$ and a Hamiltonian function $E: \Omega\rightarrow \overline{\mathbb{R}^+}$
exist, how to construct them in a natural way? }
\end{problem}

Suppose that such a canonical ensemble   $(\Omega, \mathcal{F})$ with a Hamiltonian function $E: \Omega\rightarrow \mathbb{R}^+$
exist. Let $\mathbb{P}$ be the Gibbs measure on $(\Omega, \mathcal{F})$. Let $dP(E)$ be the image measure of the Gibbs measure $\mathbb{P}$ on
$(\mathbb{R}^+, \mathcal{B}(\mathbb{R}^+))$ under $E:
\Omega\rightarrow \mathbb{R}^+$, in other words, $dP(E)$ is the
probability law of the Hamiltonian energy function $E$ considered as
a random variable from $(\Omega, \mathcal{F}, \mathbb{P})$ to $(\mathbb{R}^+, \mathcal{B}(\mathbb{R}^+))$. Then one can write
$$dP(E)={e^{-\beta E}\over Z_\beta}g(E)dE,$$
where $dE$ is the Lebesgue measure on $(\mathbb{R}^+, \mathcal{B}(\mathbb{R}^+))$, and $g(E)dE$ is the ``density of
states'' measure of the canonical ensemble  $(\Omega, \mathcal{F},
\mathbb{P})$ and the Hamiltonian $E: \Omega\rightarrow \mathbb{R}^+$. In view of this, the ``density of states'' measure
$g(E)dE$ must satisfy

\begin{eqnarray}
\log \int_{\mathbb{R}}e^{-\beta E}g(E)dE=H_n(u, t).\label{LogLap}
\end{eqnarray}This is to
say, the log-Laplace transformation of the ``density of states''
measure $g(E)dE$ on  $(\mathbb{R}^+, \mathcal{B}(\mathbb{R}^+))$
equals to $H_n(u, t)=\int_M \left({n\over 2}-f\right)dm$. Formally,
by analytic extension of the partition function $\beta\rightarrow
Z(\beta):=\int_{\mathbb{R}^+}e^{-\beta E}g(E)dE$, and using the
inverse Fourier transformation, we have (cf. \cite{Pa} p. 56)

\begin{eqnarray*}
g(E)&=&{1\over 2\pi i}\int_{\beta'-i\infty}^{\beta'+i\infty} e^{\beta E}Z(\beta) d\beta\\
&=&{1\over 2\pi}\int_{-\infty}^{+\infty} e^{(\beta'+i\beta'') E}Z(\beta'+i\beta'') d\beta'',
\end{eqnarray*}
where $\beta=\beta'+i\beta''$ s a complex variable with $\beta'>0$,
and the path of integration is parallel to the right of the
imaginary axis, i.e., along the straight line ${\rm Re}
\beta=\beta'>0$. Of course, the path may be continuously deformed so
long as the integral converges. Thanks to this observation, it might
be possible for us to prove the existence of the ``density of states
'' measure $g(E)dE$ on a natural canonical ensemble $(\Omega,
\mathcal{F}, \mathbb{P})$ and for a natural Hamitonian $E:
\Omega\rightarrow \mathbb{R}^+$ such that the partition function of
the Maxwell-Boltzmann distribution on $(\Omega, \mathcal{F}, E,
\mathbb{P})$ satisfies Perelman's condition $(\ref{LogLap})$, i.e.,
\begin{eqnarray*}
\log \int_{\mathbb{R}}e^{-\beta E}g(E)dE=H_n(u, t).
\end{eqnarray*}

\medskip

\noindent{\bf Acknowledgement}.  This paper is based on the author's
efforts for a better understanding of Perelman's mysterious
$W$-entropy for the Ricci flow. Some ideas of this paper were
obtained in June 2006 and appeared first time in the author's 2007
Habilitation Thesis submitted to the Universit\'e Paul Sabatier in
July 2007, and in the author's published papers \cite{Li11a, Li11b}.
During recent years, the author has been invited to present some related works in the 2008 Workshop of Markov Processes and
Related Topics in Wuhu (invited by Professor M.-F. Chen), the 2009
Annual Meeting of the Hong Kong Mathematical Society (invited by
Professor N. Mok), the 2010 International Conference of Stochastic
Processes and their Applications in Osaka (invited by Professors S.
Kuwae and T. Shioya), the 2010 International Congress of the Chinese
Mathematicians in Tsinghua University (invited by Professor S.-T.
Yau), the 5th International Conference on Stochastic Analysis and
its Applications in Bonn (invited by Professor K.-T. Sturm in 2011),
and the 7th Friendship Conference in Differential Geometry between
China and Japan in Tokyo (invited by Professor A. Futaki in 2012).
The author would like to thank the organizers of the above
conferences for their invitations and to thank many people for their
interests and invitations for giving seminar talks at the various
stages of this work and related works during recent years. The
author is very grateful to Professor Hua Chen and Professor Yiming
Long for inviting him to present this work and
 related works in the Nankai Conference on
Emerging Topics on Differential Equations and their Applications,
December 5-9, 2011, at the Chern Institute of Mathematics, and for
their encouragement without which it will never be possible for him
to write this paper. Part of this work was revised  when the author visited Rutgers University during April 16-22, 2012. The author would like to thank Eric Carlen for invitation, helpful discussions and comments, and to thank the colleagues of the Department of Mathematics of Rutgers University for very nice hospitality and for their interests on this work.

\medskip

\begin{flushleft}

Xiang-Dong Li\\

\medskip

{\sc  Academy of Mathematics and Systems Science, Chinese
Academy of Sciences, No. 55, Zhongguancun East Road, Beijing, 100190, China}, E-mail: xdli@amt.ac.cn
\end{flushleft}

\end{document}